\documentclass[a4paper,10pt]{llncs}

\usepackage{html}

\usepackage{graphicx}

\title{Hecke Operators for Period Functions for Congruence Subgroups}

\author{Tobias M\"uhlenbruch\thanks{This work has been supported by the Deutsche Forschungsgemeinschaft through the DFG Forschergruppe ``Zetafunktionen und lokalsymmetrische R\"aume''.}}
\institute{Institut f\"ur Theoretische Physik\\ Technische Universit\"at Clausthal \\ 38678 Clausthal-Zellerfeld \\ Germany\\ \email{tobias.muehlenbruch@tu-clausthal.de}}

\newcommand{\C}{\ensuremath{\bbbc}}
\newcommand{\Q}{\ensuremath{\bbbq}}
\newcommand{\R}{\ensuremath{\bbbr}}
\newcommand{\Z}{\ensuremath{\bbbz}}
\newcommand{\N}{\ensuremath{\bbbn}}

\newcommand{\HH}{\ensuremath{\bbbh}}
\newcommand{\HHstar}{\ensuremath{\bbbh^\ast}}
\newcommand{\RR}{\ensuremath{\mathcal{R}}}

\newcommand{\FE}{\ensuremath{\mathrm{FE}}}

\newcommand{\Gmod}{\ensuremath{{\Gamma(1)}}}

\newcommand{\Gnull}[1]{\ensuremath{{\Gamma_0\left(#1\right)}}}

\newcommand{\re}[1]{\ensuremath{{\mathrm{Re}\!\left( #1 \right)}}}
\newcommand{\im}[1]{\ensuremath{{\mathrm{Im}\!\left( #1 \right)}}}
\newcommand{\SL}[1]{\ensuremath{{\mathrm{SL}_2\!\left( #1 \right)}}}
\newcommand{\Matrix}[4]{{\textstyle \left( {#1\atop #3} \: {#2 \atop  #4} \right)}}

\begin{document}

\maketitle

\begin{abstract}
In \cite{LZ02} and \cite{DH04} period functions for the full modular group and the congruence subgroups were introduced.
It was shown that the space of period functions is in $1-1$ correspondence with the space of Maass cusp forms.

In this article we summarise the construction of Hecke operators on the period functions, based on \cite{Mu}.
In particular we compute the matrix representation of the second Hecke operator on period functions for the full modular group $\SL{\Z}$ respectively the congruence subgroup $\Gnull{2}$.
\end{abstract}

\section{Notation}
\label{A}
To state our main result and to sketch the content of each section we have to fix some notations used throughout the text.
Let $\N$ denote the set of positive integers and let $n \in \N$ throughtout this paper.
We denote by $\mathrm{M}_n(2,\Z)$ (respectively $\mathrm{M}_\ast(2,\Z)$) the set of $2\times 2$ matrices with integer entries and determinant $n$ (respectively positive determinant).
Let $\RR_n:= \Z[\mathrm{M}_n(2,\Z)]$ (respectively $\RR:= \Z[\mathrm{M}_\ast(2,\Z)]$) be the set of finite linear combinations with coefficients in $\Z$ of elements of $\mathrm{M}_n(2,\Z)$ (respectively $\mathrm{M}_\ast(2,\Z)$).
Similarly, we denote by $\mathrm{M}_n^+(2,\Z)$ (respectively $\mathrm{M}_\ast^+(2,\Z)$) the set of $2\times 2$ matrices with nonnegative integer entries and determinant $n$ (respectively positive determinant).
Denote furthermore by $\RR_n^+:= \Z[\mathrm{M}_n^+(2,\Z)]$ (respectively $\RR^+:= \Z[\mathrm{M}_\ast^+(2,\Z)]$) the set of finite linear combinations with coefficients in $\Z$ of elements of $\mathrm{M}_n^+(2,\Z)$ (respectively $\mathrm{M}_\ast^+(2,\Z)$).
Note that $\RR^{(+)} = \bigcup_{n=1}^\infty \RR_n^{(+)}$ and $\RR_n^{(+)} \cdot \RR_m^{(+)} \subset \RR_{nm}^{(+)}$.
By definition we have $\SL{\Z} = \mathrm{M}_1(2,\Z)$.
The following four elements of $\SL{\Z}$ will play a special role in our paper:
\[
I:=       \Matrix{1}{0}{0}{1}, \quad
T:=       \Matrix{1}{1}{0}{1}, \quad
S:=       \Matrix{0}{-1}{1}{0} \quad \mbox{and} \quad
T^\prime:=\Matrix{1}{0}{1}{1}.
\]
The \emph{Hecke congruence subgroup} $\Gnull{n}$ is given by
\[
\Gnull{n} :=\left\{
\Matrix{a}{b}{c}{d} \in \SL{\Z};\, c \equiv 0 \bmod n
\right\}.
\]
The full modular group is denoted by $\Gmod:=\Gnull{1}=\SL{\Z}$.
Let $\mu=\mu_n$ denote the index of $\Gnull{n}$ in $\Gmod$ and let $\alpha_1,\ldots,\alpha_\mu$ denote representatives of the right cosets in $\Gnull{n}\backslash \Gmod$.

We also need the set of upper triangular matrices
\begin{equation}
\label{A2}
X_m =
\left\{ \Matrix{a}{b}{0}{d} \in \mathrm{M}_m(2,\Z); \, d > b \geq 0 \right\}.
\end{equation}

Let $\bar{\C}$ denote the one-point compactification of $\C$.
The map
\[
\mathrm{M}_\ast(2,\Z) \times \bar{\C} \to \bar{\C};\quad
\left(\Matrix{a}{b}{c}{d},z\right) \mapsto \Matrix{a}{b}{c}{d}z :=\frac{az+b}{cz+d}.
\]
gives an action of the matrices on $\bar{\C}$.
It induces the familiar \emph{slash action} $|_s$ on functions $f$ on $\HH:=\{z\in \C;\, \im{z}>0\}$ respectively $(0,\infty)$ formally defined by
\begin{equation}
\label{A1}
f\big|_s\Matrix{a}{b}{c}{d}(z) =
(ad-bc)^s\, (cz+d)^{-2s} \,f\left(\frac{az+b}{cz+d}\right)
\end{equation}
for complex numbers $s$ and certain classes of matrices $\Matrix{a}{b}{c}{d}$.
The slash action is well defined for
\begin{description}
\item[(a)]
$s \in \Z$, $\Matrix{a}{b}{c}{d} \in \mathrm{M}_m(2,\Z)$, $m\in \N$, $z \in \HH$ and
\item[(b)]
$s \in \C$, $\Matrix{a}{b}{c}{d} \in \mathrm{M}_m^+(2,\Z)$, $m\in \N$ and $z \in (0,\infty)$
\end{description}
as the discussions in~\cite{HMM05} and \cite{Mu03} show.
We have $\big(f\big|_s\alpha\big) \big|_s \gamma = f\big|_s(\alpha\gamma)$ for all matrices $\alpha$, $\gamma \in \mathrm{M}_\ast(2,\Z)$ in case (a) respectively $\mathrm{M}_\ast^+(2,\Z)$ in case (b).
We extend the slash action linearly to formal sums of matrices.

\section{Maass cusp forms and period functions}

\begin{definition}
\label{A4}
A \emph{Maass cusp form} $u$ for the congruence subgroup $\Gnull{n}$ is a real-analytic function $u:\HH \to \C$ satisfying:
\begin{enumerate}
\item
$u(gz) =u(z)$ for all $g \in \Gnull{n}$,
\item
$\Delta u= s(1-s)u$ for some $s \in \C$ where $\Delta = -y^2(\partial_x^2+\partial_y^2)$ is the hyperbolic Laplace operator.
We call the parameter $s$ the \emph{spectral parameter} of $u$.
\item
\label{A3}
$u$ is of rapid decay in all cusps:
for $p \in \Q \cup \{\infty\}$ and $g\in\Gmod$ such that $gp=\infty$ we have $u(g^{-1}z) = \mathrm{O}\left(\im{z}^C\right)$ as $\im{z} \to \infty$ for all $C \in \R$.
\end{enumerate}
We denote the space of Maass cusp forms for $\Gnull{n}$ with spectral value $s$ by $S(n,s)$.
\end{definition}

A function $f:(0,\infty) \to \C$ is called \emph{holomorphic} if it is locally the restriction of a holomorphic function.

The vector valued period functions for $\Gnull{n}$ are defined as follows:
\begin{definition}
\label{A5}
A \emph{period function} for $\Gnull{n}$ is a function $\vec{\psi}:(0,\infty) \to \C^\mu$ with $\vec{\psi}=(\psi_i)_{i\in\{1,\ldots,\mu\}}$ such that
\begin{enumerate}
\item $\psi_i$ is holomorphic on $(0,\infty)$ for all $i \in \{1,\ldots,\mu\}$.
\item $\vec{\psi}(z) = \rho(T^{-1})\, \vec{\psi}(z+1) + (z+1)^{-2s}\,\rho({T^\prime}^{-1})\, \vec{\psi}\left(\frac{z}{z+1}\right)$.
We call the parameter $s\in \C$ the \emph{spectral parameter} of $\vec{\psi}$.
The matrix representation $\rho:\Gmod\to \C^{\mu\times \mu}$ is induced by the trivial representation of $\Gnull{n}$, see Appendix~\ref{D1}.
\item
For each $i=1,\ldots,\mu$ $\psi_i$ satisfies the growth condition
\[
\psi_i(z) = \left\{\begin{array}{ll}
\mathrm{O} \left( z^{\max\{0,-2\re{s}\}} \right) \quad& \mbox{as } z \downarrow 0 \mbox{ and}\\
\mathrm{O} \left( z^{\min\{0,-2\re{s}\}} \right) \quad& \mbox{as } z\to \infty.
\end{array} \right.
\]
\end{enumerate}
\end{definition}
Following \cite{LZ02} we denote the space of period functions for $\Gnull{n}$ with spectral value $s$ by $\FE(n,s)$.
We call a function $\vec{\psi}$ a \emph{period like function} if $\vec{\psi}$ satisfies only conditions $1$ and $2$.
The space of period like functions for $\Gnull{n}$ with spectral value $s$ is denoted by $\FE^\ast(n,s)$

\begin{remark}
The authors of \cite{LZ02} have shown that $\FE(1,s)$ is isomorphic to $S(1,s)$.
A.~Deitmar and J.~Hilgert generalize this result to submodular groups of finite index, and hence for the congruence subgroup $\Gnull{n}$, in~\cite{DH04}.
\end{remark}

Before we recall Hecke operators for $S(n,s)$ and construct Hecke operators for $\FE(n,s)$ we have to show how the spaces $S(n,s)$ and $\FE(n,s)$ are related.

\section{Vector valued cusp forms}
\label{D2}

For each $u\in S(n,s)$ we construct a vector valued version of $u$ which transforms under $\Gmod$ with respect to the representation $\rho$.

\begin{definition}
\label{D2.1}
A \emph{vector valued cusp form} $\vec{u}:\HH \to \C^\mu$ for $\Gnull{n}$ with spectral value $s \in \C$ is a vector valued function $\vec{u} = (u_1, \ldots,u_\mu)^\mathrm{tr}$ satisfying
\begin{itemize}
\item $u_j$ is real-analytic for all $j\in \{1,\ldots,\mu\}$,
\item $\vec{u}(g z) = \rho(g)\, \vec{u}(z)$ for all $z \in \HH$ and $g \in \Gmod$, where $\rho$ is the matrix representation defined in Appendix~\ref{D1},
\item $\Delta u_j = s(1-s)u_j$ for all $j\in \{1,\ldots,\mu\}$ and
\item $u_j(z) = \mathrm{O}\left(\im{z}^C\right)$ as $\im{z} \to \infty$ for all $C \in \R$ and $j \in \{1,\ldots,\mu\}$.
\end{itemize} 
We denote the space of all vector valued cusp forms with spectral parameter $s$ for $\Gnull{n}$ by $S_\mathrm{ind}(n,s)$.
\end{definition}

To each $u \in S(n,s)$ we associate the vector valued function $\Pi(u)$ given by
\begin{equation}
\label{D2.2}
\Pi: S(n,s) \to S_\mathrm{ind}(n,s); \quad u \mapsto \Pi(u):= \big( u\big|_0\alpha_1, \ldots, u\big|_0\alpha_n \big)^\mathrm{tr}.
\end{equation}
It was shown in \cite{Mu} that $\Pi$ is bijective.
This proves the following
\begin{lemma}
\label{D2.3}
The spaces $S_\mathrm{ind}(n,s)$ and $S(n,s)$ are isomorphic.
\end{lemma}

\section{The period functions of vector valued cusp forms}
\label{E3}
We identify $\pm\infty$ with the cusp $i\infty$.
The action of $\Gmod$ on $\HH$ extends naturally to $\HHstar:=\HH \cup \Q \cup \{\infty\}$.

We use the definition in \cite{Mu} of a \emph{path} $L$ connecting points $z_0,z_1 \in \HHstar$.
Basically we understand by a \emph{path} $L$ connecting points $z_0,z_1 \in \HHstar$ a piecewise smooth curve which lies inside $\HH$ except possibly for finitely many points.
If a point $z$ of the path $L$ is in $\Q \cup \{\infty\}$ then locally the path lies inside an arc bounded by geodesics starting in $z$.
For distinct $z_0,z_1 \in \HHstar \setminus \HH$ the \emph{standard path} $L_{z_0,z_1}$ is the geodesic connecting $z_0$ and $z_1$.

\begin{definition}
\label{E3.1}
For $\vec{u} \in S_\mathrm{ind}(n,s)$ and $L_{0,\infty}$ the standard path the integral transform $P:S_\mathrm{ind}(n,s) \to \big(C^\omega(0,\infty)\big)^\mu$, with $C^\omega(0,\infty)$ the space of holomorphic functions on $(0,\infty)$, is defined as
\begin{equation}
\label{E3.2}
\big(P\vec{u}\big)_j(\zeta) = \int_{L_{0,\infty}} \eta\big(u_j, R_\zeta^s \big)
\qquad \mbox{for }
\zeta >0, \, j \in \{1,\ldots,\mu\}.
\end{equation}
Formally, we write (\ref{E3.2}) as
\begin{equation}
\label{E3.2a}
\big(P\vec{u}\big)(\zeta) = \int_{L_{0,\infty}} \eta\big(\vec{u}, R_\zeta^s \big).
\end{equation}
\end{definition}

\begin{remark}
\begin{itemize}
\item
Note that $R_\zeta(z)$ denotes a function of two variables $\zeta$ and $z$.
\item
The existence of the integral transform is shown in \cite{Mu}.
\item
The $1$-forms $\eta(u_i,R_\zeta^s)$ are closed, see Appendix~\ref{E2}.
Hence 
\[
\int_{L^\prime} \eta\big(u_i, R_\zeta^s \big)
=
\int_{L_{0,\infty}} \eta\big(u_i, R_\zeta^s \big) = \big[Pu\big]_i(\zeta)
\]
for arbitrary paths $L^\prime$ homotopic to $L_{0,\infty}$.
\item
The fact that $R_\zeta^s(z)$ is holomorphic in $\zeta$, see~(\ref{E2.1}), implies that $\big[P\vec{u}\big]_j(\zeta)$ is a holomorphic function in $\zeta$.
\end{itemize}
\end{remark}

The function $P\vec{u}$ in (\ref{E3.2a}) has the following properties \cite{Mu}:

\begin{lemma}
\label{E3.3}
For $\zeta>0$ and $\gamma = \Matrix{a}{b}{c}{d} \in \Gmod$ with $a,b,c,d \geq 0$ the function $P\vec{u}$ with $\vec{u} \in S_\mathrm{ind}(n,s)$ satisfies
\begin{equation}
\label{E3.4}
(c\zeta +d)^{-2s} \, \rho(\gamma^{-1}) \, P\vec{u}(\gamma \zeta)
= 
\int_{L_{\gamma^{-1}0,\gamma^{-1}\infty}} \eta\big(\vec{u}, R_\zeta^s\big).
\end{equation}
\end{lemma}

\begin{lemma}
\label{E3.5}
For $\vec{u} \in S_\mathrm{ind}(n,s)$ each component of the function $P\vec{u}$ satisfies the growth conditions
\[
(P\vec{u})_i(\zeta) = \left\{\begin{array}{ll} 
\mathrm{O} \left( \zeta^{\max \big(0,-2\re{s}\big)} \right) \qquad & \mbox{as } \zeta \downarrow 0 \mbox{ and}\\
\mathrm{O} \left( \zeta^{\min \big(0,-2\re{s}\big)} \right) \qquad & \mbox{as } \zeta \to \infty.
\end{array}\right. 
\]
\end{lemma}

Lemma~\ref{E3.3} and Lemma~\ref{E3.5} imply the following:
\begin{proposition}
\label{E3.6}
For $\vec{u} \in S_\mathrm{ind}(n,s)$ the function $P\vec{u}$ is a period function.
\end{proposition}

Moreover, \cite{Mu} shows the following
\begin{proposition}
\label{E3.8}
For $s \in \C \setminus \Z$, $\re{s} > 0$ the operator
\[
P:\quad S_\mathrm{ind}(n,s) \to \FE(n,s); \quad \vec{u} \mapsto P\vec{u}
\]
is bijective.
\end{proposition}

\section{Hecke operators for Maass cusp forms}
\label{D5}
We define the Hecke operators on $S(n,s)$ for a fixed $s \in \C$ similarly as in~\cite{atkin:1} for modular forms.

\begin{definition}
\label{D5.2}
Denote by $T(p)$ and $U(q)$ for $\gcd(p,n)=1$, $q|n$ and $p,q$ prime the following elements in $\RR_p$ respectively $\RR_q$:
\begin{equation}
\label{D5.2a}
\label{D5.2b}
T(p) =
\sum_{ad=p \atop 0 \leq b <d} \Matrix{a}{b}{0}{d} 
\quad \mbox{and} \quad
U(q) =
\sum_{0 \leq b <q} \Matrix{1}{b}{0}{q}.
\end{equation}
The induced maps $S(n,s) \to S(n,s)$; $u \mapsto u\big|_0T(p)$ respectively $u \mapsto u\big|_0U(q)$ are called the $p^\mathrm{th}$ and $q^\mathrm{th}$ Hecke operator $H_p$ respectively $H_q$ on $S(n,s)$
\end{definition}

Obviously, the Hecke operators $H_p$ and $H_q$ depend on $n$.

\begin{remark}
The Hecke operators $H_p$ and $H_q$ use a nonstandard normalization.
\end{remark}

\section{Hecke operators for vector valued cusp forms}
\label{D4}

For each $g \in \Gmod$ we define the map $\sigma_g:X_m \to X_m$ by
$A \,g\, \left( \sigma_g(A)\right)^{-1} \in \Gmod$ for all $A \in X_m$.
It was shown in \cite{Mu} that $\sigma_g$ is bijective and satisfies $\sigma_g^{-1} = \sigma_{g^{-1}}$.

Also for $A \in X_m$ we define the map
\begin{equation}
\label{D4.8.b}
\phi_A=\phi_{A,n}: \quad \{1,\ldots,\mu_{n}\} \to \{1,\ldots,\mu_n\}; \quad
  i \mapsto \phi_A(i)
\end{equation}
such that 
\begin{equation}
\label{D4.8.c}
A \alpha_i \in \Gnull{n}\alpha_{\phi_A(i)} \,\sigma_{\alpha_i}(A).
\end{equation}
Usually we write $\phi_{A,n}=\phi_A$ omitting the index $n$ since $n$ is fixed in the entire discussion.
Note that the map $\phi_A$ is not bijective, see e.g.\ the example given in Table~\ref{ex2.1}.

\smallskip

To derive a formula for the Hecke operators acting on $S_\mathrm{ind}(n,s)$ we have to write the vector valued cusp form $\Pi\big( u\big|_0\sum_A A\big)$ in terms of a linear action of a certain matrix sum on the vector valued cusp form $\Pi(u)$.

For prime $p,q$ with $\gcd(p,n)=1$ and $q|n$ the $p^\mathrm{th}$ (respectively $q^\mathrm{th}$) Hecke operator $H_p$ (respectively $H_q$) is given by the action of $T(p)$ (respectively $U(q)$) on the space of cusp forms which we write as
\[
S(n,s) \to S(n,s); \quad
u \mapsto u \big|_0 \sum_A A = \sum_{A \in \mathcal{A}} u \big|_0 A
\]
with $\sum_A A = T(p)$ and $\mathcal{A} =X_p$ (respectively $\sum_A A = U(q)$ and $\mathcal{A}=X_q \setminus \Matrix{q}{0}{0}{1}$).

Consider the $j^\mathrm{th}$ component of the vector valued cusp form $\Pi\big( u\big|_0\sum_A A\big)$.
By (\ref{D4.8.c}) write this component as
\begin{equation}
\label{D4.10.a}
\left( u\big|_0\sum_A A\right) \big|_0\alpha_j
=
u\big|_0 \sum_A \left(A \alpha_j \right)
=
u\big|_0 \sum_A \alpha_{\phi_A(j)} \sigma_{\alpha_j}(A)
\end{equation}
for indices $\phi_A(j) \in \{1,\ldots,\mu_n\}$ and upper triangular matrices $\sigma_{\alpha_j}(A)$.

This allows us to define Hecke operators for vector valued cusp forms.
\begin{definition}
\label{D4.5}
For $n,m\in \N$, $m$ prime and $s\in \C$ put $\sum_A A := T(m)$ if $m \not\;\mid n$ and put $\sum_A A := U(m)$ if $m \mid n$.
The \emph{$m^\mathrm{th}$ Hecke operator} $H_{n,m}$ on $\vec{u} \in S_\mathrm{ind}(n,s)$ is defined as
\begin{equation}
\label{D4.5.a}
\big(H_{n,m} \vec{u}\big)_{j} \mapsto  \sum_A \, u_{\phi_A(j)} \big|_0 \sigma_{\alpha_j}(A)
\quad \mbox{for } j\in \{1,\ldots,\mu_n\}.
\end{equation}
\end{definition}
In (\ref{D4.5.a}) we sum over all $A \in X_m$ if $m \not\;\mid n$ and $A \in X_m\setminus\left\{ \Matrix{m}{0}{0}{1} \right\}$ if $m \mid n$.

\smallskip

The $m^\mathrm{th}$ Hecke operator $H_{n,m}$ on $S_\mathrm{ind}(n,s)$ corresponds to the $m^\mathrm{th}$ Hecke operator $H_m$ on $S(n,s)$ \cite{Mu}:

\begin{proposition}
\label{D4.7}
$\Pi(H_m u) = H_{n,m} \Pi(u)$.
\end{proposition}

\section{Left neighbor sequences}
\label{C2}
We recall the left neighbor sequences as introduced in \cite{Mu}.
The necessary facts on Farey sequences can be found in Appendix~\ref{C1}.

\begin{definition}
\label{C2.10}
The \emph{left neighbor map} $\mathrm{LN}: \Q \cup \{+\infty\} \to  \Q \cup \{-\infty\}$ such that $\mathrm{LN}(q)$ is the left neighbor of $q$ in the Farey sequence $F_{\mathrm{lev}(q)}$, that is
\begin{equation}
\label{C2.1}
\mathrm{LN}(q) = 
\max\{ r \in F_{\mathrm{lev}(q)}; \, r < q \}.
\end{equation}
\end{definition}

\begin{lemma}
\label{C2.2}
For $q \in \Q\cup \{+\infty\}$ and $\mathrm{lev}(q) >0$ we have $\mathrm{lev}\big(\mathrm{LN}(q)\big) < \mathrm{lev}(q)$.
\end{lemma}

\begin{definition}
\label{C2.3}
Let be $q \in \Q \cup \{+\infty\}$ and $L=L_q\in \N$ such that
\[
\mathrm{LN}^{L}(q) = -\infty
\qquad \mbox{and} \qquad
\mathrm{LN}^{l}(q) > -\infty
\quad \mbox{for all } l=1,\ldots,L-1.
\]
The \emph{left neighbor sequence} $\mathrm{LNS}(q)$ of $q$ is the finite sequence
\[
\mathrm{LNS}(q)= \big( \mathrm{LN}^{L}(q), \mathrm{LN}^{L-1}(q), \ldots, \mathrm{LN}^{1}(q), q \big),
\]
where we use the notation $\mathrm{LN}^{l}(q) := \underbrace{\mathrm{LN}\big(\mathrm{LN}(\cdots \mathrm{LN}(q)) \cdots \big)}_{l\,\mathrm{times}}$. 
\end{definition}
The number $L=L_q$ in Definition~\ref{C2.3} is unique.

To $\mathrm{LNS}(q)$ we construct an element $M(q) \in \RR_1$:
\begin{definition}
\label{C2.6}
To each rational $q \in [0,1)$ consider $\mathrm{LNS}(q)=\big(\frac{a_0}{b_0}, \ldots,\frac{a_L}{b_L} \big)$ with $\gcd(a_l,b_l)=1$ and $b_l \geq 0$, $l=0,\ldots,L$.
We define $M(q)=\sum_{l=1}^L m_l \in \RR_1$ by
\begin{equation}
\label{C2.7}
M(q) = \Matrix{-a_0}{a_1}{-b_0}{b_1}^{-1} + \ldots + 
       \Matrix{-a_{l-1}}{a_{l}}{-b_{l-1}}{b_{l}}^{-1} + \ldots +
       \Matrix{-a_{L-1}}{a_{L}}{-b_{L-1}}{b_{L}}^{-1}.
\end{equation} 
\end{definition}

\begin{lemma}
\label{C2.8}
Let $0\leq q < 1$ rational and $M(q)=\sum_{l=1}^L \Matrix{\ast}{\ast}{c_l}{d_l}$.
We have $c_l \zeta +d_l >0$ for all $\zeta > q$ and $l \in \{1,\ldots,L\}$.
\end{lemma}

\begin{lemma}
\label{C2.9}
Let $A =\Matrix{a}{b}{0}{d} \in \mathrm{M}_\ast(2,\Z)$ be such that $a,b \in \N$, $0 \leq b < d$ and $M\big(\frac{b}{d}\big)=\sum_{l=1}^L m_l$.
Then
\begin{itemize}
\item
$\det m_l =1$, 
\item
the matrices $m_lA$ contain only nonnegative integer entries and
\item
the entries of $m_lA=\Matrix{a^\prime}{b^\prime}{c^\prime}{d^\prime}$ satisfy $a^\prime > c^\prime \geq 0$ and $d^\prime > b^\prime \geq 0$.
\end{itemize}
\end{lemma}

\begin{lemma}
\label{F1.1}
For rational $q \in [0,1)$ put $M(q)=\sum_l m_l$ as in Definition~\ref{C2.6}.
The two paths $L_{q,\infty}$ and $\bigcup_l L_{m_l^{-1} 0, m_l^{-1} \infty}$ have the same initial and end point.
\end{lemma}

\section{Hecke operators for period functions for $\Gmod$}
\label{F1}
We consider first the simpler case $n=1$.

\begin{definition}
\label{F1.5}
For $m \in \N$ define
\begin{equation}
\tilde{H}(m) :=
\sum_{A \in X_m} M\big(A0\big) A \in \RR^+_m.
\end{equation}
For $s \in \C$ the formal sum $\tilde{H}(m)$ defines an operator $\tilde{H}_m$ on $C^\omega(0,\infty)$ by
\begin{equation}
\label{F1.5.b}
\tilde{H}_{1,m} f := f\big|_s \tilde{H}(m)
\qquad (f\in C^\omega(0,\infty)).
\end{equation} 
\end{definition}

\begin{remark}
\begin{itemize}
\item
Lemmas~\ref{C2.9} and \ref{C2.8} imply in particular that $\big(\tilde{H}_{1,m} f\big)(\zeta)$ is well defined for all $\zeta>0$.
\item
Lemma~\ref{C2.9} also shows that $\{m_lA; \, A \in X_m, \, M(A0)=\sum_{l=1}^L m_l\}$, containing all matrices in $\tilde{H}(m)$, is a subset of
\[
S_m=\left\{ \Matrix{a}{b}{c}{d}; \, a>c \geq 0, d>b \geq 0 \right\} \subset \mathrm{M}_n^+(2,\Z).
\]
It is shown in \cite{HMM05} that both sets are equal, implying $\tilde{H}(m) = \sum_{B\in S_m} B$.
(The authors in \cite{HMM05} assume that $\gcd(a,b,c,d)=1$ but this restriction is not necessary.\@)
\end{itemize}
\end{remark}

Before we show that $\tilde{H}_{1,m}$ are indeed the Hecke operators on $\FE(1,s)$ we need some technical lemmas \cite{Mu}:

\begin{lemma}
\label{F1.4}
Let $u \in S(1,s)$ be a cusp form and $Pu$ as in Definition~\ref{E3.1}.
For $A =\Matrix{\ast}{\ast}{0}{d} \in X_m$ and $M(A0) = \sum_{l=1}^L m_l \in \RR_1$ we have
\begin{equation}
\label{F1.4.a}
m^s d^{-2s} \, \int_{L_{A0,A\infty}} \eta\big(u, R_{A\zeta}^s\big) 
 = 
\sum_{l=1}^L \big(Pu \big|_s m_lA \big) (\zeta)= \big(Pu \big|_s M(A0)A\big)(\zeta)
\end{equation}
for all $\zeta >0$.
\end{lemma}

\begin{lemma}
\label{F1.6}
Let $Pu$ be the period function of the cusp form $u \in S(1,s)$.
For any $m \in \N$ the operator $\tilde{H}_{1,m}$ satisfies
\begin{equation}
\label{F1.6.a}
\Big(\tilde{H}_{1,m} (Pu) \Big)(\zeta)
=
\int_0^{i\infty} \eta\big(H_m u, R_\zeta^s \big)
\qquad \mbox{for } \zeta>0.
\end{equation}
\end{lemma}

\begin{proof}
Use Lemma~\ref{F1.4} together with the transformation property~(\ref{E2.2}).
\qed
\end{proof}

The relation between $H_m$ on $S(1,s)$ and $\tilde{H}_{1,m}$ on $\FE(1,s)$ is given by
\begin{proposition}
\label{F1.8}
For $u \in S(1,s)$ the period function $Pu \in \FE(1,s)$ satisfies the identity
\[
\big(Pu\big)\big|_s\tilde{H}(m) = P\big(u\big|_0 H(m)\big).
\]
\end{proposition}

\begin{proof}
This follows immediately from Lemma~\ref{F1.6}.
\qed
\end{proof}

\begin{remark}
Proposition~\ref{F1.8} shows $\tilde{H}_{1,m}$ are indeed Hecke operators on $\FE(1,s)$.
They are the same, \cite{MM05}, \cite{FMM}, as the operators constructed in \cite{HMM05} using only properties of the period functions respectively transfer operators for the groups $\Gnull{m}$.
Another derivation of $\tilde{H}_{1,m}$ for $\Gmod$ is given also in \cite{Mu04} using a criterion found by Choie and Zagier in \cite{CZ93}.
L.~Merel gives also a similar representation in \cite{Me94}.
\end{remark}

\section{Hecke operators for period functions for $\Gnull{n}$}
\label{F2}
In this section we extend the above derivation of the Hecke operators for period functions for $\Gmod$ to the congruence subgroups $\Gnull{n}$.

First, we extend Lemma~\ref{F1.4}:
\begin{lemma}
\label{F2.1}
For $A =\Matrix{\ast}{\ast}{0}{d} \in X_m$ put $M(A0) = \sum_{l=1}^L m_l \in \RR_1$.
If $P\vec{u}$ is a period function of $\vec{u} \in S_\mathrm{ind}(n,s)$, then
\begin{equation}
\label{F2.1.a}
m^s d^{-2s} \, \int_{L_{A0,A\infty}} \eta\big(\vec{u}, R_{A\zeta}^s\big)
 = 
\sum_{l=1}^L  \rho(m_l^{-1}) \, \Big( P\vec{u} \big|_s m_lA \Big) (\zeta)
\qquad \mbox{for } \zeta>0.
\end{equation}
\end{lemma}

In the following we denote the $i^\mathrm{th}$ component of the vector $\vec{u}$ by $[\vec{u}]_i$.

\begin{lemma}
\label{F2.2}
Let $\alpha_1, \ldots, \alpha_\mu$ be representatives of the right coset of $\Gnull{n}$ in $\Gmod$ where $\mu= [\Gmod:\Gnull{n}]$.
Let $P\vec{u}$ be the period function of $\vec{u} \in S_\mathrm{ind}(n,s)$.
For $A\in X_m$, $j\in \{1,\ldots,\mu_{n}\}$ let be $M(\sigma_{\alpha_j}(A)0) = \sum_{l=1}^L m_l \in \RR_1$.
Then the following identity holds for all $j \in \{1,\ldots,\mu\}$ and $\zeta >0$:
\begin{equation}
\label{F2.2.a}
\int_{L_{0,\infty}} \eta\big([\vec{u}]_{\phi_A(j)} \big|_0 \sigma_{\alpha_j}(A), R_\zeta^s \big)
=
\sum_{l=1}^L  \left[ \rho(m_l^{-1}) \, P\vec{u} \right]_{\phi_A(j)} \big|_s m_l \sigma_{\alpha_j}(A)\,(\zeta).
\end{equation}
\end{lemma}

\begin{proof}
Write $\vec{u} = (u_j)_j$.
Using Lemma~\ref{E2.5} and Property~(\ref{E2.2}) of $R_\zeta$ we find for any $j \in \{1,\ldots,\mu_n\}$ and $\zeta>0$:
\begin{eqnarray*}
&&\int_{L_{0,\infty}} \eta\big(u_{\phi_A(j)} \big|_0 \sigma_{\alpha_j}(A), R_\zeta^s \big)
\, = \,
\int_{L_{\sigma_{\alpha_j}(A)0,\sigma_{\alpha_j}(A)\infty}} \eta\big(u_{\phi_A(j)} , R_\zeta^s \big|_0 \big(\sigma_{\alpha_j}(A)\big)^{-1}\big)\\
&& \quad = \,
m^s d_j^{-2s} \, \int_{L_{\sigma_{\alpha_j}(A)0,\sigma_{\alpha_j}(A)\infty}} \eta\Big(u_{\phi_A(j)} , R_{\sigma_{\alpha_j}(A)\zeta}^s \Big)
\end{eqnarray*}
where $\sigma_{\alpha_j}(A)= \Matrix{\ast}{\ast}{0}{d_j}$ is again in $X_m$.
Take $M(\sigma_{\alpha_j}(A)0)=\sum_{l=1}^L m_l \in \RR_1$ and apply Lemma~\ref{F2.1}. 
We have
\begin{eqnarray*}
\int_{L_{0,\infty}} \eta\big( u_{\phi_A(j)} \big|_0  \sigma_{\alpha_j}(A) , R_\zeta^s \big)
&=&
\Big[ \int_{L_{0,\infty}} \eta\big( \vec{u} \big|_0  \sigma_{\alpha_j}(A) , R_\zeta^s \big) \Big]_{\phi_A(j)}\\
&=&
\left[ 
\sum_{l=1}^L  \rho(m_l^{-1}) \, \Big(P\vec{u} \Big|_s  m_l \sigma_{\alpha_j}(A) \Big)(\zeta)
\right]_{\phi_A(j)} .
\end{eqnarray*}
\qed
\end{proof}

\begin{remark}
For $u\in S(n,s)$ and $A \in X_m$ Relation~(\ref{D4.8.c}) implies that
\[
u\big|_0 A\alpha_j = u\big|_0 \alpha_{\phi_A(j)} \sigma_{\alpha_j}(A).
\]
Hence equation~(\ref{F2.2.a}) can be written as
\[
\int_{L_{0,\infty}} \eta\big(\Pi(u\big|_0 A), R_\zeta^s \big)
=
\sum_{l=1}^L  \left[ \rho(m_l^{-1}) \, P\Pi(u) \right]_{\phi_A(j)} \big|_s m_l \sigma_{\alpha_j}(A) \, (\zeta)
\]
\end{remark}

Lemma~\ref{F2.2} allows us to derive an explicit formula for the action of the Hecke operators on the period functions of $\Gnull{n}$ induced from the action of these operators on $S_\mathrm{ind}(n,s)$ for this group.

\begin{proposition}
\label{F2.3}
Let $\alpha_1, \ldots, \alpha_\mu$ be representatives of the right cosets of $\Gnull{n}$ in $\Gmod$.
Let $P\vec{u}$ be the period function of $\vec{u} \in S_\mathrm{ind}(n,s)$.
For $m$ prime take $\mathcal{A} \subset X_m$ such that $T(m)=\sum_{A \in \mathcal{A}} A$ if $\gcd(m,n)=1$ respectively $U(m)=\sum_{A \in \mathcal{A}} A$ if $m \mid n$.
The $m^\mathrm{th}$ Hecke operator $\tilde{H}_{n,m}$ acting on $P\vec{u}$ is given by
\begin{equation}
\label{F2.3.a}
\left[ \tilde{H}_{n,m} \big(P\vec{u}\big)\right]_j=
\sum_{A \in \mathcal{A}} \sum_{l=1}^L
\left[ \rho(m_l^{-1}) \, P\vec{u} \right]_{\phi_A(j)} \big|_s \big(m_l \sigma_{\alpha_j}(A) \big).
\end{equation}
\end{proposition}

\begin{remark}
We emphasize that the constant $L$ in (\ref{F2.3.a}) depends on $A$ and that $\mathcal{A} = X_m$ for $\gcd(m,n)=1$ respectively $\mathcal{A} = X_m \setminus \left\{\Matrix{m}{0}{0}{1} \right\}$ for $m \mid n$.
\end{remark}

\begin{proof}[of Proposition~\ref{F2.3}]
The $m^\mathrm{th}$ Hecke operator $H_{n,m}$ acts on $\vec{u}$ as
\[
\left[ H_{n,m} \vec{u} \right]_j =  \sum_A \, u_{\phi_A(j)} \big|_0 \sigma_{\alpha_j}(A)
\qquad \mbox{for } k \in \{1,\ldots,\mu\}.
\]
Applying Lemma~\ref{F2.2} to both sides then gives formula~(\ref{F2.3.a}).
\qed
\end{proof}

\begin{remark}
In his diploma thesis \cite{Fr}, M.~Fraczek uses a similar approach to compute a representation of the Fricke operator on $\FE(n,s)$. 
Part of his thesis was also to write a C-program computing the representation of $\tilde{H}_{n,m}$.
Moreover the approach was also used in \cite{FMM} to understand Hecke-like operators on $\FE(n,s)$ constructed in \cite{HMM05}.
\end{remark}

\begin{appendix}
\section{On Farey sequences}
\label{C1}
We recall Farey sequences and some of its properties as presented in \cite{Hu94} and \cite{Mu}.
We adhere to the convention to denote infinity in rational form as $\infty = \frac{1}{0}$ and $-\infty = \frac{-1}{0}$ and to denote rationals $\frac{p}{q}$ with coprime $p \in \Z$ and $q \in \N$.

\begin{definition}
\label{C1.1}
For $n \in \N$ the \emph{Farey sequence} $F_n$ of \emph{level $n$} is the sequence
\[
F_n := \left( \frac{u}{v};\, u,v \in \Z, |u| \leq n, 0 \leq v \leq n \right).
\]
ordered by the standard order $<$ of $\R$.
We define $F_0$ as 
\[
F_0 := \left( \frac{-1}{0}, \frac{0}{1}, \frac{1}{0} \right).
\]
The level function $\mathrm{lev}: \Q \cup \{\pm\infty\} \rightarrow \Z$ is defined by
\[
\mathrm{lev}\left(\frac{a}{b}\right)=\left\{ \begin{array}{ll} 
0 & \mbox{if } \frac{a}{b} \in \{ \frac{-1}{0}, \frac{0}{1}, \frac{1}{0} \} \mbox{ and}\\
\max\{|a|,|b|\}\qquad & \mbox{otherwise}.
\end{array}\right.
\]
\end{definition}

The number $\mathrm{lev}(q)$ is just the level of the Farey sequence in which the number $q$ appears for the first time.

\smallskip

Let $\frac{a}{c}$ and $\frac{b}{d}$ be two neighbors in the Farey sequence $F_n$.
Then the square matrix $\Matrix{a}{b}{c}{d}$ satisfies $\det \Matrix{a}{b}{c}{d} = \pm 1$.

\begin{lemma}
\label{C1.2}
Let $\frac{a}{c}$ and $\frac{b}{d}$ be two neighbors of the Farey sequence $F_n$.
Then 
\[
\det \Matrix{a}{b}{c}{d} = -1 
\quad \Longleftrightarrow \quad \frac{a}{c} < \frac{b}{d}.
\]
\end{lemma}

\begin{remark}
Our applications of the Farey sequences deal mostly with the case $\det \Matrix{a}{b}{c}{d} = -1$.
However, we prefer matrices in $\Gmod$.
For this we replace $\Matrix{a}{b}{c}{d}$ by $\Matrix{-a}{b}{-c}{d}$ which obviously does not change $\frac{a}{c}$ and $\frac{b}{d}$.
\end{remark}

\begin{lemma}
\label{C1.3}
For $\frac{a}{c}$ and $\frac{b}{d}$ with $a,b,c,d\in \Z$, $c,d \geq 0$ and $ad-bc = \pm 1$ define $n:=\max \left\{ \mathrm{lev}\left(\frac{a}{c}\right),\, \mathrm{lev}\left(\frac{b}{d}\right) \,\right\}$.
Then $\frac{a}{c}$ and $\frac{b}{d}$ are neighbors in $F_n$.
\end{lemma}

\section{Induced representations}
\label{D1}
Let $G$ be a group and $H$ be a subgroup of $G$ of finite index $\mu=[G:H]$. 
For each representation $\chi:H \to \mathrm{End}(V)$ we consider the induced representation $\chi_H:G \to \mathrm{End}(V_G)$, where
\[
V_G:=\{f:G\to V;\, f(hg)=\chi(h)f(g) \quad \mbox{for all }g\in G, h \in H\} 
\]
and
\[
\big(\chi_H(g)f\big)(g^\prime) = f(g^\prime g)
\qquad \mbox{for all }g,g^\prime \in G.
\]
For $V=\C$ and $\chi$ the trivial representation we call the induced representation $\chi_H$ the \emph{right regular representation}.
In fact, in this case $V_G$ is the space of left $H$-invariant functions on $G$ or, what is the same, functions on $H \backslash G$, and the action is by right translation in the argument. 
One can identify $V_G$ with $V^\mu$ using a set $\{\alpha_1,\ldots,\alpha_\mu\}$ of representatives for $H \backslash G$, i.e., $H \backslash G = \{H\alpha_1, \ldots, H\alpha_\mu\}$.
Then
\[
V_G \to V^\mu 
\quad \mbox{with} \quad
f \mapsto \big(f(\alpha_1),\ldots,f(\alpha_\mu)\big)
\]
is a linear isomorphism which transports $\chi_H$ to the linear $G$-action on $V^\mu$ given by
\[
g \cdot (v_1, \ldots,v_\mu) =  
\big(\chi(\alpha_1 g \alpha_{k_1}^{-1}) v_{k_1}, \ldots,  \chi(\alpha_\mu g \alpha_{k_\mu}^{-1}) v_{k_\mu}\big)
\]
where $k_j\in \{1,\ldots,\mu\}$ is the unique index such that $H\alpha_j g =H\alpha_{k_j}$. 
To see this, one simply calculates 
\[
\big( \chi_H(g)f \big)(\alpha_j) = f(\alpha_j g) = f(\alpha_j g \alpha_{k_j}^{-1} \alpha_{k_j})
= \chi(\alpha_j g \alpha_{k_j}^{-1}) \big(f(\alpha_{k_j})\big).
\]
In the case of the right regular representation the identification $V_G \cong \C^\mu$ gives a matrix realization
\[
\rho(g) = \big( \delta_H(\alpha_i g \alpha_j^{-1}) \big)_{1 \leq i,j \leq \mu}
\]
where $\delta_H(g)=1$ if $g \in H$ and $\delta_H(g)=0$ otherwise.
In particular, the matrix $\rho(g)$ is a permutation matrix.

We take $G=\Gmod$, $H=\Gnull{n}$ and $\alpha_1, \ldots, \alpha_\mu\in \Gmod$ as representatives of the $\Gnull{n}$ orbits in $\Gmod$.
The \emph{matrix representation} $\rho: \Gmod \to \C^{\mu \times \mu}$ is
\begin{equation}
\label{D1.2}
\rho(g):= 
\Big(  \delta_\Gnull{n} (\alpha_i\,g\, \alpha_j^{-1})  \Big)_{1 \leq i,j \leq \mu}
\qquad \mbox{for all }g \in \Gmod.
\end{equation}
We easily check that $\rho$ satisfies $\rho(g^\prime)\,\rho(g) = \rho(g^\prime g)$ for all $g, g^\prime \in \Gmod$.

\section{The function $R_\zeta(z)$ and the 1-form $\eta(\cdot,\cdot)$}
\label{E2}
We define the function $R_\zeta(z):\C \times \HH \to \C$ as
\begin{equation}
\label{E2.1}
R_\zeta(z) = \frac{y}{(x-\zeta)^2+y^2}.
\end{equation}
Note that $R_\zeta(z)$ is a nonstandard notation for a function of the two variables $\zeta \in \C$ and $z\in \HH$.

For $\zeta \in \R$ we have $R_\zeta(z) = \frac{\im{z}}{|z-\zeta|^2}$.

It is shown in \cite{LZ02}, \cite{Mu03} and also in \cite{Mu} that $R_\zeta(z)$ satisfies the transformation formula
\begin{equation}
\label{E2.2}
\frac{|\det g|}{|c\zeta+d|^2 } \,R_{g\zeta}(gz) = R_\zeta(z)
\end{equation} 
for all $g \in \mathrm{M}_1(2,\Z)$ and real $\zeta$.
Moreover, it was also shown there that $R_\zeta^s(z)$ is an eigenfunction of the hyperbolic Laplace operator with eigenvalue $s(1-s)$ for all complex $s$.
% , \cite{LZ02},
% \begin{equation}
% \label{E2.3}
% \Delta \big(R_\zeta(z)\big)^s =  s(1-s)\, \big(R_\zeta(z)\big)^s
% \qquad (s \in \C).
% \end{equation}

\smallskip

As in \cite{LZ02}, we define the $1$-form $\eta(u,v)$ for two smooth functions $u,v$ on $\HH$:
\begin{equation}
\label{E2.4}
\eta(u,v) := \big(v \partial_yu-u\partial_yv \big)dx + \big(u \partial_xv-v\partial_xu \big)dy.
\end{equation} 
The following lemma is shown in \cite{LZ02}:
\begin{lemma}
\label{E2.5}
If $u$ and $v$ are eigenfunctions of $\Delta$ with the same eigenvalue, then the $1$-form $\eta(u,v)$ is closed.
If $ z \mapsto g(z)$ is any holomorphic change of variables, then the $1$-form satisfies $\eta(u \circ g,v \circ g) =\eta(u,v) \circ g$.
\end{lemma}

\section{Example: $\tilde{H}_{1,2}$}
\label{ex1}
For suitable complex $s$ and $u \in S(1,s)$ let $Pu$ be the period function of $u$ given in \S\ref{E3}.
Consider the $2^\mathrm{nd}$ Hecke operator $H_2$ on $S(1,s)$ given by
\[
u \mapsto H_2 u = u \big|_0 T(2).
\]
In this example we compute the matrix representation $\tilde{H}(2)$ of the $2^\mathrm{nd}$ Hecke operator $\tilde{H}_{1,2}$ on $\FE(1,s)$.

Recall that the set $X_2$, see (\ref{A2}), and the element $T(2)$, see (\ref{D5.2a}), are given by
\begin{equation}
X_2 = \left\{\Matrix{1}{0}{0}{2}, \Matrix{1}{1}{0}{2}, \Matrix{2}{0}{0}{1}\right\}
\; \mbox{ and } \;
T(2) = \sum_{A\in X_2}A = \Matrix{1}{0}{0}{2} + \Matrix{1}{1}{0}{2} + \Matrix{2}{0}{0}{1}.
\end{equation}

Given the Farey sequences $F_0$, $F_1$ and $F_2$ we compute $\mathrm{LNS}(A0)$ and the $M(A0)$ for all $A \in X_2$, see \S\ref{C2}:
\begin{eqnarray}
\label{ex1.5}
&&
\mathrm{LNS}(0)=\left( \frac{-1}{0},\frac{0}{1}\right), \quad 
\mathrm{LNS}\big({\textstyle \frac{1}{2}}\big)=\left( \frac{-1}{0},\frac{0}{1},\frac{1}{2}\right),\\
\label{ex1.6}
&&M(0)=\Matrix{1}{0}{0}{1} 
\quad \mbox{and} \quad
M\big({\textstyle \frac{1}{2}}\big)=\Matrix{1}{0}{0}{1} + \Matrix{2}{-1}{1}{0}.
\end{eqnarray}

For $P\left(u\big|_0A\right)(\zeta)$, $\zeta>0$, with $A=\Matrix{1}{1}{0}{2} \in X_2$ we find
\[
P\left(u\big|_0A\right)(\zeta)
=
\int_{L_{0,\infty}} \eta(u\big|_0A ,R_\zeta^s) (z)
=
\int_{L_{A0,A\infty}} \eta(u ,R_\zeta^s\big|_0A^{-1}) (z)
\]
where the slash operator acts on the $z$ variable of $R_\zeta(z)$.
The transformation property~(\ref{E2.2}) gives
\[
P\left(u\big|_0A\right)(\zeta)
=
2^s \, \int_{L_{A0,A\infty}} \eta\big(u ,R_{A\zeta}^s\big) (z)
=
2^s \, \int_{L_{A0,A\infty}} \eta\left(u ,R_\frac{\zeta+1}{2}^s\right) (z).
\]
We would like to write $\int_{L_{A0,A\infty}} \eta(u ,R_{A\zeta}^s) (z) = \sum_l Pu(m_lA\zeta)$.
This can be done using $\mathrm{LNS}\big({\textstyle \frac{1}{2}}\big)$ and the related formal sum $M\big({\textstyle \frac{1}{2}}\big)$.
We have
\[
L_{A0,A\infty} = L_{\frac{1}{2},\infty} 
=
L_{0,\infty}
\cup
L_{\frac{1}{2},0}
= 
L_{m_1^{-1}0,m_1^{-1}\infty}
\cup
L_{m_2^{-1}0,m_2^{-1}\infty}
\]
with $m_1=I \in \Gmod$ and $m_2=\Matrix{2}{-1}{1}{0}\in \Gmod$.
Hence, $P\left(u\big|_0A\right)(\zeta)$ can be written as
\begin{eqnarray*}
P\left(u\big|_0A\right)(\zeta)
&=&
2^s \, \int_{L_{A0,A\infty}} \eta(u ,R_{A\zeta}^s) (z)\\
&=&
2^s \, \int_{L_{m_1^{-1}A0,m_1^{-1}A\infty}} \eta(u ,R_{m_1A\zeta}^s) (z)\\
&& +
2^s (\zeta+1)^{-2s} \, \int_{L_{m_2^{-1}A0,m_2^{-1}A\infty}} \eta(u ,R_{m_2A\zeta}^s) (z)\\
&=&
(Pu)\big|_s (m_1A+m_2A) (\zeta) = (Pu)\big|_s \big(M(A0)A\big)
\end{eqnarray*} 
where $M\big(\frac{1}{2}\big)=m_1+m_2$ with $m_1=I$ and $m_2=\Matrix{2}{-1}{1}{0}$ as in (\ref{ex1.6}).

We perform the analogous computation for the other matrices in $X_2$ and find that
\begin{equation}
\label{B1.3}
P\big(u\big|_0T(2)\big) = \big(Pu\big)\big|_s \tilde{H}(2)
\end{equation}
where the formal matrix sum $\tilde{H}(2)$ has the form
\begin{eqnarray}
\label{ex1.2}
\tilde{H}(2)
&=& 
\sum_{A\in X_2} M(A0)A
=
M(0) \, \Matrix{1}{0}{0}{2} + M \big({\textstyle \frac{1}{2}}\big) \, \Matrix{1}{1}{0}{2} 
+ M(0) \, \Matrix{2}{0}{0}{1}\\
\nonumber
&=& 
I \,\Matrix{1}{0}{0}{2} + \left[  \Matrix{1}{0}{0}{1} + \Matrix{2}{-1}{1}{0}\right] \Matrix{1}{1}{0}{2} 
+ I\, \Matrix{2}{0}{0}{1}\\
\nonumber
&=&
\Matrix{1}{0}{0}{2} + \Matrix{1}{1}{0}{2} + \Matrix{2}{0}{1}{1}  + \Matrix{2}{0}{0}{1}.
\end{eqnarray}

\section{Example: $\tilde{H}_{2,2}$}
\label{ex2}
For suitable complex $s$ and $u \in S(2,s)$ let $P\Pi u$ be the period function of $u$ given in \S\ref{E3}.
Consider the $2^\mathrm{nd}$ Hecke operator $H_2$ on $S(2,s)$ given by
\[
u \mapsto H_2 u = u \big|_0 U(2).
\]
In this example, we compute the matrix representation $\tilde{H}(2)$ of the $2^\mathrm{nd}$ Hecke operator $\tilde{H}_{2,2}$ on $\FE(2,s)$.

According to (\ref{D5.2b}) we have
\[
U(2) = \Matrix{1}{0}{0}{2} + \Matrix{1}{1}{0}{2} = A_1 + A_2 \in \RR_2
\]
using $A_1=\Matrix{1}{0}{0}{2}$, $A_2=\Matrix{1}{1}{0}{2}$ and $A_3=\Matrix{2}{0}{0}{1}$.

We already computed $M(A_10)=M(A_30)=M(0)$ and $M(A_20)=M\big(\frac{1}{2}\big)$ in (\ref{ex1.5}).

Put
\begin{equation}
\label{ex2.5}
\alpha_1=I=\Matrix{1}{0}{0}{1}, \quad
\alpha_2=S=\Matrix{0}{-1}{1}{0} \quad \mbox{and} \quad
\alpha_3=ST=\Matrix{0}{-1}{1}{1}
\end{equation}
as representatives of the right cosets in $\Gnull{2}\backslash \Gmod$.
The values of the functions $\phi_A=\phi_{A,2}$ and $\sigma_{\alpha_i}$ defined in \S\ref{D4} are given in Table~\ref{ex2.1}.

The matrix representation $\rho$, see (\ref{D1.2}), has the values
\begin{equation}
\label{C2.5}
\rho(I)=\left(\begin{array}{ccc}
1&0&0\\
0&1&0\\
0&0&1
\end{array}\right) 
\qquad \mbox{and} \quad
\rho\left( \Matrix{2}{-1}{1}{0}^{-1} \right)=\left(\begin{array}{ccc}
0&1&0\\
1&0&0\\
0&0&1
\end{array}\right).
\end{equation}

\smallskip

Computing $\tilde{H}_{2,2}$ via (\ref{F2.3.a}), we find for the first component
\begin{eqnarray}
\nonumber
\left[ \tilde{H}_{2,2}\big(P\vec{u}\big)\right]_1 &=&
\sum_{A\in \{A_1, A_2\}} \sum_{l=1}^{L_{\sigma_{\alpha_1}(A)}}
\left[ \rho(m_l^{-1}) \, P\vec{u} \right]_{\phi_A(1)} \big|_s \big(m_l \sigma_{\alpha_1}(A) \big)\\
\nonumber
&=& 
\left[ \rho(I) \, P\vec{u} \right]_{\phi_{A_1}(1)} \big|_s \big(I \sigma_{\alpha_1}(A_1) \big)
+
\left[ \rho(I) \, P\vec{u} \right]_{\phi_{A_2}(1)} \big|_s \big(I \sigma_{\alpha_1}(A_2) \big)\\
\nonumber
&&\quad
+
\left[ \rho\left(\Matrix{2}{-1}{1}{0}^{-1}\right) \, P\vec{u} \right]_{\phi_{A_2}(1)} \big|_s \big(\Matrix{2}{-1}{1}{0} \sigma_{\alpha_1}(A_2) \big)\\
\nonumber
&=&
\left[ P\vec{u}\right]_1 \big|_s IA_1 +
\left[ P\vec{u}\right]_1 \big|_s IA_2 + \left[ P\vec{u}\right]_2 \big|_s \Matrix{2}{-1}{1}{0} A_2 \\
\label{ex2.3}
&=&
\left[ P\vec{u}\right]_1 \big|_s \Matrix{1}{0}{0}{2} + 
\left[ P\vec{u}\right]_1 \big|_s \Matrix{1}{1}{0}{2} + \left[ P\vec{u}\right]_2 \big|_s \Matrix{2}{0}{1}{1}.
\end{eqnarray}
Similarly, we find for the second and third component
\begin{eqnarray}
\label{ex2.4}
\left[ \tilde{H}_{2,2}\big(P\vec{u}\big)\right]_2 
& = &
\left[ P\vec{u}\right]_2 \big|_s \Matrix{2}{0}{0}{1} + 
\left[ P\vec{u}\right]_1 \big|_s \Matrix{1}{1}{0}{2} + \left[ P\vec{u}\right]_2 \big|_s \Matrix{2}{0}{1}{1}
\mbox{ and}\\
\label{C2.4a}
\left[ \tilde{H}_{2,2}\big(P\vec{u}\big)\right]_3 
& = &
\left[ P\vec{u}\right]_2 \big|_s \Matrix{2}{0}{0}{1} + 
\left[ P\vec{u}\right]_1 \big|_s  \Matrix{1}{0}{0}{2} .
\end{eqnarray}

\begin{remark}
\begin{itemize}
\item 
Note that the third component of $\tilde{H}_{2,2}\big(P\vec{u}\big)$ has only two terms compared to the three terms of the other components.
This is related to the fact that $A_2$ does not occur in the image of $\sigma_{\alpha_3}(\{A_1,A_2\})$.
\item
The operator $\tilde{H}_{2,2}$ does not use the third component of the vector valued period form since the index $3$ does not appear in the image of $\phi_{A}$ in Table~\ref{ex2.1}.
\item
M.~Fraczek, \cite{Fr}, wrote a C-program for the computation of $\tilde{H}_{n,m}$.
\end{itemize} 
\end{remark}

\begin{table}
\begin{center}
\begin{tabular}{c||c|c||c|c||c|c}
$j$& $\phi_{A_1}(j)$ &$\sigma_{\alpha_j}(A_1)$ & $\phi_{A_2}(j)$ &$\sigma_{\alpha_j}(A_2)$ & $\phi_{A_3}(j)$ &$\sigma_{\alpha_j}(A_3)$\\
\hline
$1$ & $1$ & $A_1$ & $1$ & $A_2$ & $1$ & $A_3$\\ 
$2$ & $2$ & $A_3$ & $1$ & $A_2$ & $2$ & $A_1$\\ 
$3$ & $2$ & $A_3$ & $1$ & $A_1$ & $2$ & $A_2$\\ 
\end{tabular}
\caption{The values of the functions $\phi_A=\phi_{A,2}$ and $\sigma_{\alpha_j}(A)$ for $m=n=2$ and $A \in X_2$. 
The representatives $\alpha_j$ of the cosets of $\Gnull{2}\backslash \Gmod$ are given in~(\ref{ex2.5}).}
\label{ex2.1}
\end{center}
\end{table}

\end{appendix}


\begin{thebibliography}{1}

\bibitem{atkin:1}
A.~O.~L. Atkin and J.~Lehner.
Hecke operators on $\Gamma_0(m)$.
Math.\ Ann.\ \textbf{185} (1970) 134--160.

\bibitem{CZ93}
Y.~J.~Choie and D.~Zagier.
Rational period functions for $\mathrm{PSL}(2,\bbbz)$.
In M.~Knopp and M.~Sheingorn, editors, {\em A Tribune to Emil Grosswald: Number Theory and related Analysis}, Volume 143 of {\em  Contemporary Mathematics}, pages 89--108.
American Mathematical Society, 1993.

\bibitem{DH04}
A.~Deitmar and J.~Hilgert.
The Lewis Correspondence for submodular groups.
e-arxiv (2004).
\texttt{http://arXiv.org/abs/math/0404067}.

\bibitem{Fr}
M.~Fraczek.
Spezielle Eigenfunktionen des Transfer-Operators f\"ur Hecke Kongruenz Untergruppen.
Diploma thesis (2006), TU Clausthal.

\bibitem{FMM}
M.~Fraczek, D.~Mayer, T.~M\"uhlenbruch.
A realization of the Hecke algebra on the space of period functions for $\Gamma_0(n)$.
J.~Reine Angew.~Math. accepted.
\texttt{http://de.arxiv.org/abs/math.NT/0512355}

\bibitem{HMM05}
J.~Hilgert, D.~Mayer and H.~Movasati.
Transfer operators for {$\Gamma_0(n)$} and the Hecke operators for period functions of $\mathrm{PSL}(2,\bbbz)$.
Math.~Proc.~Camb.~Phil.~Soc. \textbf{139} (2005) 81--116.

\bibitem{Hu94}
A.~Hurwitz.
Ueber die angen\"aherte Darstellung der Zahlen durch rational Br\"uche.
Mathematische Annalen \textbf{44} (1894) 417--436.

% \bibitem{La76}
% S.~Lang.
% \emph{Introduction to modular forms}.
% Volume 222 of Grundl.\ math.\ Wiss., Springer Verlag, Berlin, 1976.

\bibitem{LZ02}
J.~Lewis and D.~Zagier.
Period functions for Maass wave forms. {I}.
Ann.\ of Math.\ \textbf{153} (2001) 191--258.

\bibitem{MM05}
D.~Mayer, T.~M\"uhlenbruch.
From the transfer operator for geodesic flows on modular surfaces to the Hecke operators on period functions of $\Gamma_0(n)$.
In \emph{Dynamical Systems: from Algebraic to Topological Dynamics}.
Proceedings of the ESF-Exploratory Workshop, 5-9 July 2004, Bonn.
Contemporary Mathematics \textbf{385} (2005), 137--161, American Mathematical Society.

\bibitem{Me94}
L.~Merel.
Universal Fourier expansions of modular forms.
In G.~Frey, editor, {\em On Artin's conjecture for odd $2$-dimensional representations}, volume 1585 of Lecture Notes in Math.
Springer Verlag, Berlin, 1994.

\bibitem{Mi89}
T. Miyake.
\emph{Modular Forms}.
Springer-Verlag, 1989.

\bibitem{Mu03}
T.~M\"uhlenbruch.
Systems of automorphic forms and period functions.
PhD thesis, Utrecht University, September 2003.

\bibitem{Mu04}
T.~M\"uhlenbruch.
Hecke operators on period functions for the full modular group.
IMRN 4127--4145 (2004).

\bibitem{Mu}
T.~M\"uhlenbruch.
Hecke operators on period functions for $\Gamma_0(N)$.
J.\ Number Theory. accepted.
\texttt{http://dx.doi.org/10.1016/j.jnt.2005.09.003}

\end{thebibliography}
\end{document}